\documentclass[11pt]{amsart}

\usepackage[centertags]{amsmath}
\usepackage{graphicx}

\numberwithin{equation}{section}

\theoremstyle{plain}
\newtheorem{theorem}{Theorem}[section]
\newtheorem{corollary}[theorem]{Corollary}
\newtheorem{lemma}[theorem]{Lemma}
\newtheorem{proposition}[theorem]{Proposition}

\theoremstyle{definition}

\newtheorem{definition}[theorem]{Definition}

\theoremstyle{remark}
\newtheorem{remark}[theorem]{Remark}

\newcommand{\J}{\mathbb{J}}
\newcommand{\R}{\mathbb{R}}

\newcommand{\N}{\mathbb{N}}
\newcommand{\PP}{\mathbb{P}}

\newcommand{\CS}{\mathcal{S}}

\begin{document}

\title[Random mappings]{Asymptotic evolution of acyclic random mappings}

\author{Steven N.\ Evans}
\address{Steven N.\ Evans \\
  Department of Statistics \#3860 \\
  University of California at Berkeley \\
  367 Evans Hall \\
  Berkeley, CA 94720-3860 \\
  U.S.A.}
\email{evans@stat.Berkeley.EDU}
\urladdr{http://www.stat.berkeley.edu/users/evans/}
\thanks{SNE supported in part by NSF grant DMS-0405778}

\author{Tye Lidman}
\address{Tye Lidman}
\email{tlid@berkeley.edu}

\date{\today}



\begin{abstract} 
An acyclic mapping from an $n$ element set into itself
is a mapping $\varphi$ such that if $\varphi^k(x) = x$
for some $k$ and $x$, then $\varphi(x) = x$.  Equivalently,
$\varphi^\ell = \varphi^{\ell+1} = \ldots$ for $\ell$ sufficiently large.
We investigate the behavior as 
$n \rightarrow \infty$ of a Markov chain
on the collection of such mappings. At each step of the chain,
a point in the $n$ element set is chosen uniformly at random
and the current mapping is modified by replacing
the current image of that point by a new one
chosen independently and uniformly at random, conditional
on the resulting mapping being again acyclic. 
We can represent an acyclic mapping as a directed graph
(such a graph will be a collection of rooted trees)
and think of these directed graphs as metric
spaces with some extra structure. Heuristic
calculations indicate that the metric space valued 
process associated with the Markov chain should,
after an appropriate time and ``space'' rescaling,
converge as $n \rightarrow \infty$ to a real tree ($\R$-tree)
valued Markov process that is reversible with respect
to a measure induced naturally by the standard reflected
Brownian bridge.  The limit process, which we construct
using Dirichlet form methods, is a Hunt
process with respect to a suitable Gromov-Hausdorff-like
metric. This process is similar to
one that appears in earlier work by Evans and Winter
 as the limit of chains involving the
subtree prune and regraft tree (SPR) rearrangements from
phylogenetics.
\end{abstract}
\maketitle

\section{Introduction}
\label{S:intro}

\begin{figure}[htbp]
	\begin{center}
		\includegraphics[width=1.00\textwidth]{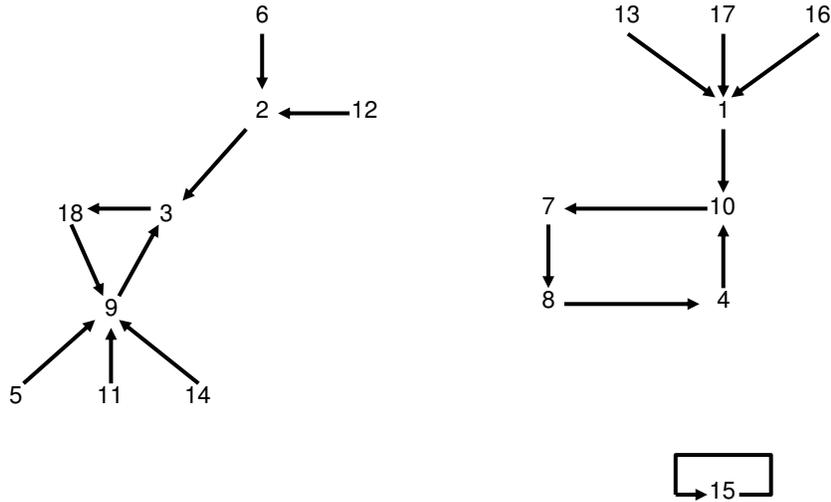}
	\end{center}
	\caption{The directed graph corresponding to the mapping in Table~\ref{table:random_map_ex}.}
	\label{fig:random_map_ex}
\end{figure}

A mapping $\varphi$ from the set $[n] := \{1,2,\ldots,n\}$ into itself may be represented as
a directed graph with vertex set $[n]$ and directed edges of the form $(i,\varphi(i))$, $i \in [n]$.
The resulting directed graph has the feature that every vertex has out-degree $1$ (with self-loops -- corresponding to fixed points
-- allowed), and any such graph corresponds to a unique mapping.  For example, the  mapping
$\varphi:[18] \rightarrow [18]$ in Table~\ref{table:random_map_ex} corresponds to the directed graph in 
Figure~\ref{fig:random_map_ex}.

\begin{table}[htbp]
{\small
\begin{center}
		\begin{tabular}{c|c|c|c|c|c|c|c|c|c|c|c|c|c|c|c|c|c|c}
		$i$   & 1 & 2 & 3 & 4 & 5 & 6 & 7 & 8 & 9 & 10 & 11 & 12 & 13 & 14 & 15 & 16 & 17 & 18 \\
		$\varphi(i)$ & 10& 3 & 18& 10& 9 & 2 & 8 & 4 & 3 & 7  & 9  & 2  & 1  & 9  & 15 & 1  &  1 & 9 \\
		\end{tabular}
	\end{center}
}
	\caption{A mapping from $[18]$ into itself.}
	\label{table:random_map_ex}
\end{table}

The directed graph may be decomposed into a number of connected
components.  Each of these
components consists of a single directed cycle (possibly a self-loop)
plus trees rooted
at each vertex on the directed cycle (such a tree may be a
trivial tree consisting of only the root, meaning that the only
pre-image of that point is its predecessor on the directed cycle).  
We call such rooted trees the {\em tree components}
of the graph.

\begin{figure}[htbp]
	\begin{center}
		\includegraphics[width=1.00\textwidth]{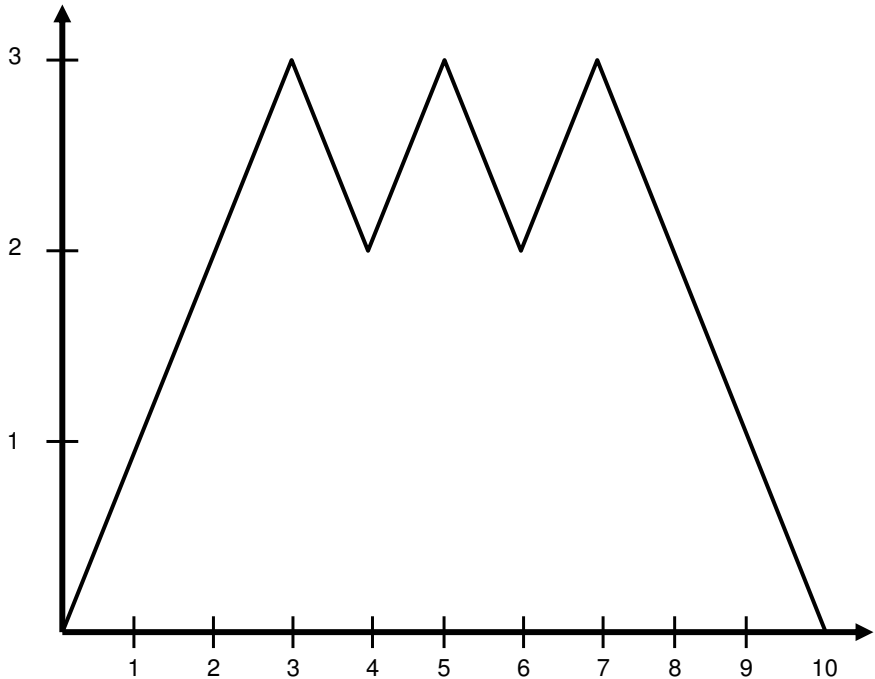}
	\end{center}
	\caption{The excursion corresponding to the tree component rooted at vertex $10$
	in Figure~\ref{fig:tree_component_walk} with the start of the excursion
	shifted to time $0$.}
	\label{fig:tree_component_walk}
\end{figure}

Aldous and Pitman \cite{MR1293075} describe a procedure for associating a mapping
of $[n]$ into itself with a {\em lattice reflected bridge path} of length $2n$, that is, with
a function $b:\{0,1,\ldots,2n\} \rightarrow \{0,1,2,\ldots\}$ such that
$b(0) = b(2n) = 0$ and $|b(k+1) - b(k)| = 1$ for $0 \le k < 2n$.  The exact details
of the procedure aren't important for us.  However, we note that  a tree component with $\ell$ vertices
corresponds to 
a lattice positive excursion path from $0$
with $2 \ell$ steps.  Such a segment of path records the distance from the root plus $1$ in a depth-first-search
of the tree component. 
For example, the tree component of
 size $5$ consisting of the vertices
$\{1,10,13,16,17\}$ in Figure~\ref{fig:random_map_ex} corresponds 
to the excursion shown in Figure~\ref{fig:tree_component_walk} after a suitable
translation of the time axis.
In particular, a tree component that consists of just one point (which is necessarily
a point on a directed cycle) 
corresponds to an excursion of the form
$b(k-1) = 0$, $b(k) = 1$, and $b(k+1) = 0$.

Of course, the mapping cannot be recovered from just the lattice reflected bridge path.
For one thing, some extra marking of distinguished
points of the zero set of the lattice path is required to split the
lattice path up into sub-paths corresponding to
components of the directed graph.  Once this is done, the 
mapping is uniquely specified by the lattice
path up to a relabeling of the vertices: that is, if two mappings
$\varphi$ and $\psi$ correspond to the same lattice reflected bridge path, then
$\psi = \pi \circ \varphi \circ \pi^{-1}$ for some permutation $\pi$ of
$[n]$.  Conversely, if $\psi = \pi \circ \varphi \circ \pi^{-1}$
for some permutation $\pi$ of $[n]$, then the lattice path
corresponding to $\psi$ may be obtained from the lattice path
corresponding to $\varphi$ by composition with
a bijective map of $\{0,1,\ldots, 2n\}$ that preserves lengths
of excursions above all levels. That is, if the lattice path corresponding
to $\varphi$ has $k$ excursions above some level $h$, then the same is
true of the lattice path corresponding to $\psi$.

Suppose now that a mapping of $[n]$ into itself is chosen
uniformly at random from the $n^n$ possibilities. This is equivalent
to choosing the image of each point of $[n]$ independently and uniformly
at random from $[n]$.  The corresponding lattice reflected bridge path is
not uniformly distributed.  However, it is shown in \cite{MR1293075} that
if the lattice reflected bridge path is turned into
a continuous time process by holding it
constant between integer time points,  
time is rescaled by $2n$, 
and space is rescaled by $n^{\frac{1}{2}}$ to  produce a 
function from $[0,1]$ into $\R_+$, then this stochastic process
with c\`adl\`ag  sample paths converges in distribution to twice
a standard reflected Brownian bridge (that is, twice the Brownian
bridge reflected at $0$ that goes from position
$0$ at time $0$ to position $0$ again at time $1$).
In particular, the proportion of vertices that lie on directed cycles converges
to the proportion of time the standard reflected Brownian bridge
spends at $0$, which is, of course, $0$, so that asymptotically
almost all vertices are not roots of tree components.  
The asymptotics of the cyclic vertices jointly with the
tree vertices are described in
\cite{MR1293075} using the local time at $0$ of the reflected bridge and that paper
also describes an auxiliary ``marking'' procedure for describing
the joint asymptotics of the the component sizes. Some later results
in this same vein may be found in
\cite{MR2197134, MR1913728, MR1445035, MR1700003,
MR1768499, MR2081482, MR1877634}.

A mapping $\varphi$ from $[n]$ into itself is {\em acyclic} if the only directed cycles in
the corresponding directed graph are self-loops.  That is, each $x \in [n]$ is either 
a fixed point of $\varphi$ (so that $x$ is a vertex on a self-loop)
or $\varphi^k(x) \ne x$ for any $k$.  Equivalently, 
$\varphi^\ell = \varphi^{\ell+1} = \ldots$ for $\ell$ sufficiently large.  For such a mapping,
each graph component consists of single tree component with
a self-loop attached to the root, and no auxiliary marking procedure
is necessary to recover the mapping up to a permutation from the corresponding lattice reflected bridge path.  It is not
hard to show that if we turn the
lattice reflected bridge path for a uniformly chosen acyclic random mapping into a
continuous time process indexed by $[0,1]$ as above, then the resulting process
also converges to twice a standard reflected Brownian bridge
-- as one would expect from the observation that the cyclic vertices are
asymptotically negligible for a uniformly chosen random mapping, 

In this paper we are interested in the asymptotic behavior 
as $n \rightarrow \infty$ of a simple
Markov chain that randomly evolves an acyclic mapping from $[n]$ into itself.
At each step of the chain, a point of $[n]$ is chosen uniformly
at random and the image of this point is re-set to a new image chosen
independently and uniformly at random from $[n]$, conditional
on the resulting mapping being acyclic.  It is clear for
each  $n$ that
this chain is reversible with respect to the uniform distribution
on the set of acyclic mappings from $[n]$ into itself and that
the chain converges to this distribution at large times.

In terms of the corresponding directed graphs, the chain
evolves as follows.  A directed edge is first chosen uniformly at random
and deleted. The deleted edge is then replaced by another directed edge with
the same initial vertex but a uniformly chosen final vertex, conditional
on the resulting graph having no cycles other than self-loops.  
Note that
the effect of such a step is the following.
\begin{itemize}
\item
If the deleted edge is a self-loop, its deletion
turns the graph component that contained the
edge into a rooted subtree.  Otherwise,
the deletion of the directed edge splits the graph component that contained
it into two pieces, one of which contains a self-loop and the other of
which is a rooted subtree.
\item
In either case, the addition of the new directed edge either attaches the
root of the subtree 
to itself by a self-loop, producing an extra graph component, or the new directed edge
 attaches the root  to a vertex chosen uniformly outside the subtree
(possibly to a vertex outside the subtree but within the same former graph component).
All such possibilities are equally likely. 
\end{itemize}

The effect on the corresponding lattice
bridge path is to remove an excursion above some level,
insert a suitable time-space translation of it at some time
point in the lattice bridge path outside the excursion, and then
close up the gap left by the removal (more precisely, this transformation
may need to be followed by a bijective map of $\{0,1,\ldots, 2n\}$ 
that preserves lengths of excursions above all levels because
of the way that the labeling of vertices in the directed graph
is used to construct the corresponding lattice bridge path).

In order to understand the asymptotic behavior of this 
sequence of chains
as $n \rightarrow \infty$, we need to embed the state space
of each chain into a common state space that will also
be the state space of the limit process.

To begin with, we erase all of the self-loops in the directed graph corresponding
to an acyclic mapping.  This produces a forest of subtrees rooted at vertices
that were formerly on self-loops.
We then connect
the roots of these subtrees by directed edges to a single adjoined point to produce a tree
rooted at the adjoined point.
Keeping in mind the rescaling
identified by Aldous and Pitman, we think of this
tree as a one-dimensional cell complex by regarding each
edge as a segment of length $n^{-\frac{1}{2}}$.  We thus
have a metric space with a distinguished base point (the
root).  This pointed metric space is an instance of
a {\em rooted compact real tree ($\R$-tree)}: 
see Section~\ref{S:realtrees} for the
precise definition of a $\R$-tree -- for the moment,
all that is important for explaining
our results is that a $\R$-tree is a metric space
that is, in some sense, ``tree-like''.
We regard two  rooted compact $\R$-trees as being equal
if one can be mapped into the other by an isometry
that preserves the root.  If two mappings $\varphi$ and $\psi$
are related by a relabeling $\psi = \pi \circ \varphi \circ \pi^{-1}$
for some permutation $\pi$ of $[n]$, then they
correspond to the same rooted compact $\R$-tree.

Before we continue with the motivation of our results,
we need to indicate how the rooted compact $\R$-tree associated with
a mapping from $[n]$ into itself may be constructed 
directly from the corresponding
lattice reflected bridge path.  We begin by introducing
some general notation that is useful later.

\begin{definition}
\label{D:Omega+}
Write $C(\R_+, \R_+)$ for the space of 
continuous functions from $\R_+$ into $\R_+$.
For $f \in C(\R_+, \R_+)$, put 
\[
\zeta(f) := \inf\{s>0 : f(t) = 0 \; \text{for all $t>s$}\}
\]
with the usual convention that $\inf \emptyset = \infty$.
The set of {\em positive bridge paths} is the set 
$\Omega_+ \subset C(\R_+, \R_+)$ given by
\[
   \Omega_+
 :=
   \left\{
    f\in C(\R_+, \R_+):\,
   \begin{array}{cc} 
   f(0)=0,\, 0< \zeta(f) < \infty,\\
   \text{$f(t) \ge 0$  for $0<t<\zeta(f)$}.
   \end{array}
   \right\}
\] 
For $\ell > 0$, 
set $\Omega_+^\ell := \{f \in \Omega_+ : \zeta(f) =\ell\}$.
\end{definition}\smallskip

We associate each $f \in \Omega_+^1$ with a compact metric space as follows.
Define an equivalence relation $\sim_f$ on $[0,1]$ by letting
\[
   u_1\sim_f u_2,\,\quad\mbox{ iff }\quad\,f(u_1)=\inf_{u\in
     [u_1\wedge u_2,u_1\vee u_2]}f(u)=f(u_2).
\]
Consider the
pseudo-metric $d_{T_f}$ on $[0,1]$ defined by
\[
   d_{T_f}(u_1,u_2)
 :=
   f(u_1)-2\inf_{u\in
   [u_1\wedge u_2,u_1\vee u_2]}f(u)+f(u_2).
\]
This pseudo-metric becomes a true metric on the quotient space 
$T_f:= [0,1] /_{\sim_f}$.
The resulting metric space is compact and is an instance
of a rooted compact $\R$-tree if we define the root to be the image of 
$0$ under the quotient map.

Suppose that the function $f \in \Omega_+^1$ is obtained by first
linearly interpolating
the lattice reflected bridge path associated with an
acyclic mapping $\varphi$ of
$[n]$ into itself to produce a function in $\Omega_+^{2n}$ and then rescaling
time by $2n$ and space by $n^{\frac{1}{2}}$.  The corresponding pointed metric
space $T_f$ is the rooted compact $\R$-tree associated with 
$\varphi$ that we described above.

Any metric space of the form $T_f$ for $f \in \Omega_+^1$ has
two natural Borel measure on it.  Firstly, there is the ``uniform''
probability measure $\nu_{T_f}$ given by the push-forward
of Lebesgue measure on $[0,1]$ by the quotient map.  We call
this measure the {\em weight} on $T_f$.  Secondly,
there is the natural {\em length} measure 
$\mu_{T_f}$, which is the one-dimensional
Hausdorff measure associated with the metric $d_{T_f}$ restricted
to points of $T_f$ that are not ``leaves'' (see Section~\ref{S:realtrees} for a more precise 
definition).  When $f$ is associated with a map of $[n]$ into itself
as above, then $\mu_{T_f}$ is just the ``Lebesgue measure'' on the cell
complex $T_f$ that assigns mass $n^{-\frac{1}{2}}$ to each edge of $T_f$
(recall that we have rescaled so that each edge has length $n^{-\frac{1}{2}}$).

Now, if we speed up time by a factor of $n^{\frac{1}{2}}$ in our
Markov chain for evolving mappings of $[n]$ into itself and look
at the corresponding rooted compact $\R$-tree-valued process, then
it is reasonable at the heuristic level that we should obtain in the
limit as $n \rightarrow \infty$ a continuous time
Markov process with the following informal description.
The state space of the limit process is the space consisting
of rooted compact $\R$-trees $T$ equipped with a probability measure
$\nu_T$: we call such objects weighted rooted compact $\R$-trees. We note
that, as in the special case of $\R$-trees of the form $T_f$ for $f \in \Omega_+^1$,
an arbitrary compact $\R$-tree has a canonical length measure $\mu_T$
given by the restriction 
of the one-dimensional
Hausdorff measure associated with the metric
to the set of points that aren't leaves.  The process evolves
away from its state at time $0$ by choosing
a point $(t,v)$ at rate $dt \otimes \mu_T(dv)$ in time and
on the current tree $T$, and at time $t$
the subtree above $v$ (that is, the subtree
of points on the other side of $v$ from the root) is re-attached at a point
$w$  chosen according to $\nu_T$ (conditional on $w$ being outside the
subtree).

In general, the measure $\mu_T$ may have infinite total mass.  For example,
if $f \in \Omega_+^1$ is chosen according to the distribution of
standard reflected Brownian bridge, so that
$T_{2 f}$ is the rooted compact $\R$-tree that arises
from a limit as $n \rightarrow \infty$ of uniform acyclic random mappings
of $[n]$ into itself, then $\mu_{T_{2f}}$ almost surely
has infinite total mass.  Consequently, the above specification
of the dynamics of the limit process does not make rigorous
sense for general weighted rooted compact
$\R$-trees.  The aim of this paper is to use Dirichlet form methods
to construct a suitably well-behaved Markov process with evolution
dynamics that conform to the heuristic description.  

We do not, however,
obtain a convergence result.  The limit process has no obvious Feller-like properties
and it is not clear how to define its dynamics for all starting points (as opposed
to almost all starting points with respect to the symmetrizing
measure, which is all the Dirichlet form approach provides) in
such a way that, say, martingale problem methods might be used
to establish convergence.

The process we construct is somewhat similar to the process constructed in \cite{MR2243874}
as a limit a natural chain based on the {\em subtree prune and regraft 
(SPR)} tree rearrangement transformations from phylogenetics.
Both processes involve the relocation of a subtree whose root is chosen
according to the length measure on the current tree.  However,
the state space of the process in \cite{MR2243874} consists of weighted
unrooted compact $\R$-trees, whereas we work with weighted rooted compact $\R$-trees and
the root plays a crucial role in defining the dynamics.  The symmetrizing
measures are, as a consequence, rather different: the measure in
\cite{MR2243874} is the distribution of the Brownian continuum random tree, which is the
$\R$-tree ``inside'' twice a standard Brownian excursion, whereas
our symmetrizing measure is the distribution of the 
$\R$-tree ``inside'' twice a standard reflected Brownian bridge.
However, many of the steps in the construction are quite similar
so we omit several arguments and simply refer to the analogous ones in \cite{MR2243874}.

We note that Markov processes with reflected bridge paths as their state space
and continuous sample paths have been studied in
\cite{MR1959795, MR1891060, MR1814427}.  These processes
are reversible with respect to the distribution of a Bessel bridge 
of some index.

\section{Weighted $\R$-trees}
\label{S:realtrees}

\begin{definition}
A metric space $(T,d)$  is a {\em real tree} ($\R$-tree)
if it satisfies the following
axioms.

\noindent{\bf Axiom~0:} The space $(T,d)$ is complete.

\noindent{\bf Axiom~1:} For all 
$x,y\in T$ there exists a unique isometric embedding
  $\phi_{x,y}:[0,d(x,y)]\to T$ such that $\phi_{x,y}(0)=x$
and $\phi_{x,y}(d(x,y))=y$.

\noindent{\bf Axiom~2:} For every injective continuous map 
$\psi:[0,1]\to T$ one
has $\psi([0,1])=\phi_{\psi(0),\psi(1)}([0,d(\psi(0),\psi(1))])$.
\end{definition}\smallskip

Axiom~1 says simply  that there is a unique ``unit speed'' path between
any two points $x$ and $y$.  We write $[x,y]$ for the
image of this path and call it the {\em segment} with
{\em endpoints} $x$ and $y$.  Axiom~2 implies that the image
of any injective path connecting two points $x$ and $y$ coincides with
the segment $[x,y]$, and so such a path
may be re-parameterized
to become the unit speed path.  Thus, while Axiom~1
is satisfied by many other spaces such as $\R^d$
with the usual metric, Axiom~2 captures the
essence of ``treeness'' and is only satisfied
by $\R^d$ when $d=1$.  See
\cite{Dre84, DreTer96, MR1399749, MR1379369, Ter97, MR1851337} for 
background on $\R$-trees.  In particular, \cite{MR1851337}
shows that a number of other definitions are equivalent
to the one above.  See also \cite{Eva06}, where much of this
material is synthesized and combined with other material on probability
on $\R$-trees.

\begin{figure}[htbp]
	\begin{center}
		\includegraphics[width=1.00\textwidth]{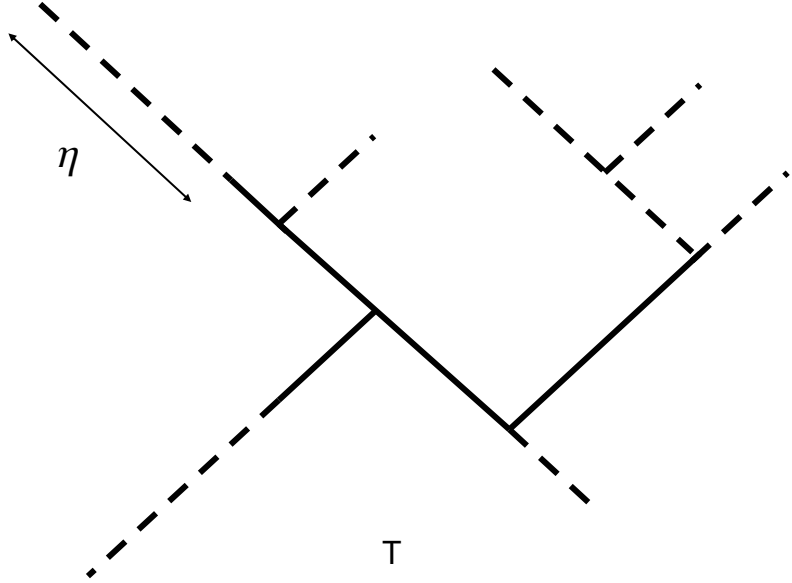}
	\end{center}
	\caption{A $\R$-tree $T$ and its $\eta$-trimming $R_\eta(T)$.
	The $\R$-tree $T$ consists of both the solid and dashed segments, whereas
	the $\R$-tree $R_\eta(T)$ consists of just the solid segments.}
	\label{fig:trimming}
\end{figure}

We define the {\em $\eta$-trimming}, $R_\eta(T)$ of a compact 
$\R$-tree $(T,d)$ for $\eta>0$ to be the set of
points $x \in T$ such that $x$ belongs to a segment
$[y,z]$ with $d(x,y) = d(x,z) = \eta$ -- see Figure~\ref{fig:trimming}.  
The {\em skeleton} of $(T,d)$ is the set
$T^o := \bigcup_{\eta>0} R_\eta(T)$.  Thus $x \in T^o$
if $x \in ]y,z[$ for some $y,z$.  The {\em leaf set}
of $(T,d)$ is the set $T \backslash T^o$.
The {\em length measure}
on $T$ is the $\sigma$-finite measure $\mu_T$ 
on the Borel $\sigma$-field ${\mathcal B}(T)$
given by the trace onto $T^o$ of the one-dimensional
Hausdorff measure associated with $d$.  Equivalently, $\mu_T$ is the unique measure
concentrated on $T^o$ such that $\mu_T([x,y]) = d(x,y)$ for all $x,y \in T$
(see Section 2.4
of \cite{MR2221786} or Section 2 of \cite{MR2243874}).

In the following, we are interested in
compact $\R$-trees $(T,d)$ equipped with a distinguished base
point $\rho \in T$ (called the {\em root}) and 
a probability measure
$\nu$ on the Borel $\sigma$-field ${\mathcal B}(T)$ (called the {\em weight}).
We call such objects {\em weighted rooted compact $\R$-trees}.
We say that two weighted rooted compact $\R$-trees $(X,d_X,\rho_X,\nu_X)$ and 
$(Y,d_Y,\rho_Y,\nu_Y)$ are {\em weighted rooted isometric} if there 
exists a bijective isometry $\Phi$ between the metric spaces $(X,d_X)$ and $(Y,d_Y)$ 
such that $\Phi(\rho_X) = \rho_Y$
and the {\em push-forward} of $\nu_X$ by $\Phi$ is $\nu_Y$, that is,
\[
   \nu_Y = \Phi_\ast\nu_X := \nu_X\circ\Phi^{-1}.
\]
The property of being weighted rooted isometric
is an equivalence relation.  
We write ${\mathbf T}^{\mathrm{wr}}$ for the collection of
equivalence classes of weighted rooted compact $\R$-trees.

In order to define a metric on
${\mathbf T}^{\mathrm{wr}}$, we first
recall the definition of the Prohorov distance between
two probability measures (see, for example, 
\cite{MR838085}).  Given two probability measures $\alpha$
and $\beta$ on a metric space $(X,d)$ with the  corresponding
collection of closed sets denoted by ${\mathcal{C}}$,
the Prohorov distance between
them is
\[
d_{\mathrm{P}}(\alpha,\beta)
:= 
\inf\{\varepsilon>0: \alpha(C) \le \beta(C^\varepsilon) + \varepsilon 
\text{ for all } C \in {\mathcal{C}}\},
\]
where $C^\varepsilon := \{x \in X: \inf_{y \in C} d(x,y) < \varepsilon\}$.
The Prohorov distance is a metric on the collection of probability
measures on $X$.

We are now in a position to define the weighted rooted Gromov-Hausdorff distance 
between the two  weighted rooted compact 
$\R$-trees $(X,d_X,\rho_X,\nu_X)$ and $(Y,d_Y,\rho_Y,\nu_Y)$.  

For  $\varepsilon>0$,
let $F^\varepsilon_{X,Y}$ denote the set of Borel maps $f:X \rightarrow Y$
such that $f(\rho_X) = \rho_Y$ and 
\[
\sup\{|d_X(x',x'') - d_Y(f(x'),f(x''))|: x', x'' \in X\} \le \varepsilon,
\]
and define $F^\varepsilon_{Y,X}$ similarly.
Put
\[
\begin{aligned}
 &  \Delta_{\mathrm{GH^{wr}}}(X,Y)
  \\
 &:=
   \inf\Bigg\{\varepsilon>0:\,\begin{array}{cc}\mbox{exist }
   f\in F^\varepsilon_{X,Y},g\in F^\varepsilon_{Y,X}\mbox{ such that} \\
   d_{\mathrm{P}}(f_\ast\nu_X,\nu_Y)\le
   \varepsilon,\,d_{\mathrm{P}}(\nu_X,g_\ast\nu_Y)\le
   \varepsilon\end{array}\Bigg\}.
\end{aligned}
\]
Note that the set on the right hand
side is non-empty because $X$ and $Y$ are compact, and hence
bounded in their respective metrics.  Note also that 
$\Delta_{\mathrm{GH^{wr}}}(X,Y)$ only depends on the weighted rooted
isometry classes of $X$ and $Y$.

It turns out that the function $\Delta_{\mathrm{GH^{wr}}}$
satisfies all the properties of a metric except the triangle inequality.
To rectify this, put
\[
d_{\mathrm{GH^{wr}}}(X,Y)
:= \inf\left\{\sum_{i=1}^{n-1} \Delta_{\mathrm{GH^{wr}}}(Z_i,Z_{i+1})^{\frac{1}{4}}\right\},
\]
where the infimum is taken over all finite sequences of weighted rooted compact
$\R$-trees $Z_1, \ldots Z_n$ with $Z_1 = X$ and $Z_n = Y$ (the exponent
$\frac{1}{4}$ is not particularly important, any sufficiently small
number would suffice).  
Note again that 
$d_{\mathrm{GH^{wr}}}(X,Y)$ only depends on the weighted rooted
isometry classes of $X$ and $Y$.

From now on, we think of $\Delta_{\mathrm{GH^{wr}}}$
and $d_{\mathrm{GH^{wr}}}$ as being defined on 
${\mathbf T}^{\mathrm{wr}} \times {\mathbf T}^{\mathrm{wr}}$.
Parts (i) and (ii) of the following result are analogous
to Lemma 2.3 of \cite{MR2243874}, part (iv) is analogous
to Proposition 2.4 of \cite{MR2243874}, part (v)
is a re-statement of Lemma 2.6 of \cite{MR2243874}, and  part (vi)
is analogous to Theorem 2.5 of \cite{MR2243874}.  
The results in \cite{MR2243874} are for $\R$-trees
with weights but without roots, but the addition
of roots does not present any new difficulties
({\em cf.} the passage from $\R$-trees without
weights or roots to $\R$-trees without weights but
with roots in Section 2.3 of \cite{MR2221786}).
The space ${\mathbf T}$ in part (iv) is the collection
of isometry classes of compact $\R$-trees (without
weights or roots) and we
refer the reader to Section 2.1 of \cite{MR2221786}
for the definition of the associated Gromov-Hausdorff
distance $d_{\mathrm{GH}}$. \smallskip

\begin{proposition}
\label{P:metr}
\begin{itemize}
\item[(i)]
The map 
$\Delta_{\mathrm{GH^{wr}}}$ has the properties:
\begin{itemize}
\item[(a)]
$\Delta_{\mathrm{GH^{wr}}}(X,Y) = 0$ if and only if $X=Y$,
\item[(b)]
$\Delta_{\mathrm{GH^{wr}}}(X,Y) = \Delta_{\mathrm{GH^{wr}}}(Y,X)$.
\end{itemize}
\item[(ii)] 
The map 
$d_{\mathrm{GH^{wr}}}$ is a metric on ${\mathbf T}^{\mathrm{wr}}$.  
\item[(iii)] For all $X,Y \in {\mathbf T}^{\mathrm{wr}}$,
\[
\frac{1}{2} \Delta_{\mathrm{GH^{wr}}}(X,Y)^{\frac{1}{4}}
\le
d_{\mathrm{GH^{wr}}}(X,Y) 
\le \Delta_{\mathrm{GH^{wr}}}(X,Y)^{\frac{1}{4}}.
\]
\item[(iv)]
A subset $D$ of $({\mathbf T}^{\mathrm{wr}},d_{\mathrm{GH^{wr}}})$ 
is relatively compact if and only if the subset
$ E :=\{(T,d):\,(T,d,\rho,\nu)\in D \}$ of
$({\mathbf T},d_{\mathrm{GH}})$ is relatively compact.
\item[(v)]
A subset $ E $ of $({\mathbf T},d_{\mathrm{GH}})$
is relatively compact if and only if
\[
\sup\{\mu_T(R_\eta(T)) : T \in  E \} < \infty
\]
for all $\eta > 0$.
\item[(vi)]
The metric space
$({\mathbf T}^{\mathrm{wr}},d_{\mathrm{GH^{wr}}})$ is 
  complete and separable.
\end{itemize}
\end{proposition}\smallskip

We note that an extensive study of spaces of metric spaces
equipped with measures is given in \cite{MR2237206, MR2237207},
and the theory of weak convergence for random variables
taking values in such spaces is developed in \cite{math.PR/0609801}.

\section{Trees and continuous paths}
\label{S:treesandpaths}

\begin{definition}
\label{D:Omega++}
The space of {\em positive excursion paths} is the set 
$\Omega_{++} \subset \Omega_+ \subset C(\R_+, \R_+)$ given by
\[
   \Omega_{++}
 :=
   \left\{
    f\in C(\R_+, \R_+):\,
   \begin{array}{cc} 
   f(0)=0,\, 0 < \zeta(f) < \infty,\\
   \text{$f(t) > 0$  for $0<t<\zeta(f)$}.
   \end{array}
   \right\}
\] 
For $\ell > 0$, 
set $\Omega_{++}^\ell := \{f \in \Omega_{++} : \zeta(f) =\ell\}$.
\end{definition}\smallskip

The following result is a slight generalization of Lemma 3.1 in \cite{MR2243874}.
The latter result was for the special case of $\R$-trees
constructed from positive excursion paths rather than
general positive bridge paths.  The proof goes through unchanged.

\begin{lemma}\label{e:comT} 
For each $f\in \Omega_+^1$, the metric space
  $(T_f,d_{T_f})$
is a compact $\R$-tree. 
\end{lemma}\smallskip

We root a $\R$-tree $(T_f,d_{T_f})$ coming from a positive bridge path in $f \in \Omega_+^1$
by taking the root to be the point corresponding to $0 \in [0,1]$ under
the quotient map.   
We equip $(T_f,d_{T_f})$ with the weight
$\nu_{T_f}$ given by the push-forward of Lebesgue measure on $[0,1]$ by the quotient map.

For a positive bridge path 
$f \in \Omega_+^1$, we identify the length measure $\mu_{T_f}$ on the
associated compact $\R$-tree $(T_f,d_{T_f})$ as follows (the discussion
is essentially the same as that in Section~3 of \cite{MR2243874}
which considered $\R$-trees coming from positive excursion paths).
For $a \ge 0$, let
\begin{equation}
\label{defGa}
   {\mathcal G}(f,a)
 := 
   \left\{s \in [0,1]:\,\begin{array}{cc}
   \text{$f(s) = a$ and, for some $t>s$,}\\ 
   \text{$f(r)>a$ for all $r \in ]s,t[$,}\\
   f(t)=a.\end{array}\right\}
\end{equation}
denote the countable set of starting points of excursions of the
function $f$ above the level $a$.  
Then, the length measure $\mu_{T_f}$ is the
push-forward of the measure
\begin{equation}
\label{length_on_path}
m_f := \int_0^\infty \mathrm{d}a \, \sum_{t \in {\mathcal G}(f,a)} \delta_t
\end{equation}
by the quotient map, where $\delta_t$ is the unit point mass at $t$.  

Alternatively, write 
\[
   \Gamma(f):=\{(s,a):\; s\in]0,1[,\,a\in[0,f(s)[\}
\]
for the region between the time axis and the graph of $f$, and for
$(s,a) \in \Gamma(f)$ denote by 
\begin{equation}
\label{start}
\underline{s}(f,s,a) := \sup\{r < s : f(r) = a\}
\end{equation}
and
\begin{equation}
\label{finish}
\bar s(f,s,a) := \inf\{t >  s : f(t) = a\}
\end{equation}
the start and finish of the
excursion of $e$ above level $a$ that straddles time $s$.  Then,
\begin{equation}
\label{length_on_path_alt}
m_f = \int_{\Gamma(f)} \mathrm{d}s \otimes \mathrm{d}a \, \frac{1}{\bar s(f,s,a) - \underline{s}(f,s,a)} 
\delta_{\underline{s}(f,s,a)}.
\end{equation}

\section{A path transformation connecting reflected
Brownian bridge and Brownian
excursion}

Write $\PP_+$ for the law of the standard Brownian bridge
reflected at $0$ that
goes from $0$ at time $0$ to $0$ at time $1$. Write
$\PP_{++}$ for the law of standard Brownian excursion.
Thus, $\PP_+$ is a probability measure on $\Omega_+^1$ and
$\PP_{++}$ is a probability measure on $\Omega_{++}^1$.
We show in this section how various computations for 
$\PP_+$ can be reduced to computations for $\PP_{++}$
using a result of Bertoin and Pitman.

Given $f \in \Omega_+^\ell$, put
\[
L(t;f) := 
\begin{cases}
\limsup_{\varepsilon \downarrow 0} \frac{1}{2 \varepsilon}
\int_{[0,t]} ds \, \mathbf{1} \{f(s) < \varepsilon\}, 
& \text{if } 
\limsup_{\varepsilon \downarrow 0} \frac{1}{2 \varepsilon}
\int_{[0,\ell]} ds \,  \mathbf{1} \{f(s) < \varepsilon\} \\
& < \infty, \\
0, 
& \text{otherwise,}
\end{cases}
\] 
for $0 \le t \le \ell$, and set $L(t;f) = L(\ell;f)$ for $t \ge \ell$.

Denote by $\tilde \Omega_+^\ell$ the subset of $\Omega_+^\ell$
consisting of functions $f$ with the properties:
\begin{itemize}
\item
the closed set $\{t \in [0,\ell] : f(t) = 0\}$ is perfect (that is,
has no isolated points) and has Lebesgue measure zero;
\item
for $0 \le t \le \ell$,
\[
L(t;f) 
=
\lim_{\varepsilon \downarrow 0} \frac{1}{2 \varepsilon}
\int_{[0,t]} ds \,  \mathbf{1} \{f(s) < \varepsilon\};
\]
\item
the function $t \mapsto L(t;f)$ is continuous; 
\item
the set of points of increase of the function $t \mapsto L(t;f)$
coincides with $\{t \in [0,\ell] : f(t) = 0\}$.
\end{itemize}
Note that if $f \in \tilde \Omega_+^\ell$, then $L(\cdot;f)$
is not identically $0$ (indeed, $L(\cdot;f)$ has $0$ as a point
of increase).  Of course, $\PP_+(\tilde \Omega_+^1) = 1$.

For $f \in \Omega_+^\ell$, set
\[
U(f) := \sup\left\{0 \le t \le \ell: L(t;f) \le \frac{1}{2} L(\ell;f)\right\}
\]
and put
\[
K^\rightarrow(t;f) :=
\begin{cases} 
L(t;f),& 0 \le t \le U(f), \\
L(\ell;f) - L(t;f),&  U(f) \le t \le \ell, \\
0,& t \ge \ell.
\end{cases}
\]
For $f \in \Omega_+^\ell$ and $u \in [0,\ell]$, set
\[
K^\leftarrow(t;f,u)
:=
\begin{cases} 
\min_{t \le s \le u} f(s),& 0 \le t \le u, \\
\min_{u \le s \le t} f(s),& u \le t \le \ell, \\
0,& t \ge \ell.
\end{cases}
\]

The following result is elementary and we leave the proof to the reader.

\begin{lemma}
\label{L:K_arrow_relations}
Fix a function $f \in \tilde \Omega_+^\ell$.
Set 
\[
e = K^\rightarrow(\cdot;f) + f.
\]
Then, $e \in \Omega_{++}^\ell$ and 
\[
f = K^\leftarrow(\cdot; e, U(f)).
\]
\end{lemma}\smallskip

The next result, which is Lemma 3.3 of \cite{MR1268525}, says that
under $\PP_+$ the path-valued random variable
$f \mapsto K^\rightarrow(\cdot;f) + f$ has law $\PP_{++}$, 
the random variable $f \mapsto U(f)$ is uniformly distributed on $[0,1]$,
and these two random variables are independent.

\begin{proposition}
\label{P:Bertoin_Pitman}
For any Borel function 
$F: \Omega_+^1 \times [0,1] \rightarrow \R_+$,
\[
\int \PP_+(\mathrm{d}f) \, F(K^\rightarrow(\cdot;f) + f, U(f))
=
\int \PP_{++}(\mathrm{d}e)  \int_{[0,1]} \mathrm{d}u \, \, F(e,u).
\]
\end{proposition}\smallskip

In order to apply Proposition~\ref{P:Bertoin_Pitman}, we need to 
understand for a fixed positive bridge path $f \in \tilde \Omega_+^1$ how the measure $m_f$
of (\ref{length_on_path}) or (\ref{length_on_path_alt}) is related
to the analogous measure for the associated positive excursion path
$K^\rightarrow(\cdot;f) + f \in \Omega_{++}^1$.

\begin{definition}
For $e \in \Omega_{++}^1$ and $u \in [0,1]$, write
\[
\Gamma^*(e,u) := \{(s,a) \in \Gamma(e) : 
u \notin [\underline{s}(e,s,a), \bar{s}(e,s,a)]\}
\]
for the set of points in $\Gamma(e)$ such that the
corresponding straddling sub-excursion does not
straddle the time $u$.
\end{definition}\smallskip

\begin{lemma} 
\label{L:ident_length}
Fix $f \in \tilde \Omega_+^1$. Set $e = K^\rightarrow(\cdot;f) + f \in \Omega_{++}^1$,
so that 
\[
\Gamma^*(e, U(f)) = 
\{(s,a) : s\in]0,1[, \, 
K^\rightarrow(s;f) \le a < K^\rightarrow(s;f) + f(s)\}.
\]
Define a bijection $\xi: \Gamma(f) \rightarrow \Gamma^*(e, U(f))$
by setting 
\[
\xi(s,a) := (s, a + K^\rightarrow(s;f)).
\]
The map $\xi$ is a measure-preserving bijection
between the set $\Gamma(f)$ equipped with the measure 
\[
\mathrm{d}s \otimes \mathrm{d}a \, \frac{1}{\bar s(f,s,a) - \underline{s}(f,s,a)}
\]
and the set $\Gamma^*(e,U(f))$ equipped with the measure
\[
\mathrm{d}s \otimes \mathrm{d}a \, \frac{1}{\bar s(e,s,a) - \underline{s}(e,s,a)}.
\]
\end{lemma}\smallskip

\begin{proof}
Decompose the open set $\{t \in [0,1] : f(t) > 0\}$
into a countable union of intervals $A_k$, $k \in \N$.
Set $B_k = \{(s,a) \in \Gamma(f) : s \in A_k\}$, $k \in \N$,
and $C_k = \{(s,a) \in \Gamma^*(e,U(f)) : s \in A_k\}$, $k \in \N$.
We have $\lambda([0,1] \setminus \bigcup_k A_k) = 0$, where $\lambda$
is Lebesgue measure.
Thus,
$\lambda \otimes \lambda(\Gamma(f) \setminus \bigcup_k B_k) = 0$
and
$\lambda \otimes \lambda(\Gamma^*(e,U(f)) \setminus \bigcup_k C_k) = 0$.

The function $t \mapsto L(t;f)$ is constant on each of the sets $A_k$,
and so the same is true of the function $t \mapsto K^\rightarrow(t;f)$.
Write $c_k$ for this constant. The function $\xi$ maps
$B_k$ bijectively into $C_k$ and the restriction of $\xi$
to $B_k$ is the translation $(s,a) \mapsto (s, a + c_k)$.

Therefore, $\xi$ is a measure-preserving bijection
between the set $\Gamma(f)$ equipped with the measure 
$\mathrm{d}s \otimes \mathrm{d}a$
and the set $\Gamma^*(e,U(f))$ equipped with the measure
$\mathrm{d}s \otimes \mathrm{d}a$.

It remains to note that if, for some $(s,a) \in \Gamma(f)$, we write $\xi(s,a) = (s,a')$, 
then we have
$\underline{s}(f,s,a) = \underline{s}(e,s,a')$
and
$\bar{s}(f,s,a) = \bar{s}(e,s,a')$, so that, in particular,
$\bar{s}(f,s,a) - \underline{s}(f,s,a) 
= \bar{s}(e,s,a') - \underline{s}(e,s,a')$.
\end{proof}\smallskip

\begin{remark}
Assume that 
$f \in \tilde \Omega_+^1$. For $a \ge 0$, recall the definition of
${\mathcal G}(f,a)$ from (\ref{defGa}).
For $u \in [0,1]$ and $e \in \Omega_{++}^1$ put
\[
   {\mathcal G}^*(e,a,u)
 := 
   \left\{s \in [0,1]:\,\begin{array}{cc}
   \text{$e(s) = a$ and, for some $t>s$,}\\ 
   \text{$e(r)>a$ for all $r \in ]s,t[$,}\\
   e(t)=a,\\
   u \notin [s,t].\end{array}\right\}
\]
That is, ${\mathcal G}^*(e,a,u)$ is the countable set of starting points
of excursions of $e$ above the level $a$ that don't straddle the time $u$.
A consequence of Lemma~\ref{L:ident_length} is that the measure
$m_f$ coincides with the measure
\[
\int_0^\infty \mathrm{d}a \, \sum_{t \in {\mathcal G}^*(K^\rightarrow(\cdot;f) + f, a, U(f))} \delta_t.
\]
\end{remark}\smallskip

As explained in the Introduction, the dynamics of the process we wish to construct
involves ``picking'' a point $v$ in a rooted compact $\R$-tree $(T,d_T,\rho_T)$ according to the
length measure $\mu_T$ and then re-rooting the subtree above $v$ (that is,
the subtree consisting of points of $x \in T$ such that $v \in [\rho_T,x[$) at a new location
$w$.  When $T = T_f$ for some $f \in \Omega_+^1$, this re-rooting of a subtree
corresponds to a rearrangement of $f$ by relocating an excursion of $f$ above some level.
We introduce the following notation to describe such rearrangements.

\begin{definition}
For $f \in \Omega_+^1$ and $(s,a) \in \Gamma(f)$, define
$\hat f^{s,a} \in \Omega_{++}$ and $\check f^{s,a} \in \Omega_{+}$, 
by
\[
   \hat f^{s,a}(t) 
 :=
\begin{cases}
   f(\underline{s}(f,s,a) + t) - a,& 0 \le t \le \bar s(f,s,a) - 
   \underline{s}(f,s,a), \\
   0,& t > \bar s(f,s,a) - \underline{s}(f,s,a),
\end{cases}
\]
and
\[
   \check f^{s,a}(t)
 :=
\begin{cases}
f(t),&  0 \le t \le \underline{s}(f,s,a),\\
f(t+ \bar s(f,s,a) - \underline{s}(f,s,a)),&  t > \underline{s}(f,s,a).
\end{cases}
\]
That is, $\hat f^{s,a}$ is the sub-excursion of $f$ that straddles $(s,a)$
shifted to start at position $0$ at time $0$, and  $\check f^{s,a}$
is $f$ with the sub-excursion that straddles $(s,a)$ excised and the
resulting gap closed up.
\end{definition}\smallskip

\begin{definition}
For $f \in \Omega_+^1$, $u \in [0,1]$, and $(s,a) \in \Gamma^*(f,u)$,
put
\[
\check U(f,u,s,a) = 
\begin{cases}
u,& 0 \le u < \underline{s}(f,s,a),\\
u - \bar s(f,s,a) + \underline{s}(f,s,a),& \bar s(f,s,a) < u \le 1.
\end{cases}
\]
By definition of $\Gamma^*(f,u)$, the point $u$ belongs to the set 
\[
[0,\underline{s}(f,s,a)[ \, \cup \, ]\bar{s}(f,s,a),1]
\]
of length 
\[
\zeta(\check f^{s,a}) = 1 - \zeta(\hat f^{s,a}) = 1 - (\bar{s}(f,s,a) - \underline{s}(f,s,a)),
\]
and $\check U(f,u,s,a)$ is where $u$ is moved to when we close up the gap
to form the interval $[0,\zeta(\check f^{s,a})]$.
\end{definition}\smallskip

The following result is immediate from Lemma~\ref{L:K_arrow_relations} and Lemma~\ref{L:ident_length}.

\begin{corollary}
\label{C:length_bridge_length_excursion}
Fix $f \in \tilde \Omega_+^1$. Set $e = K^\rightarrow(\cdot;f) + f \in \Omega_{++}^1$.
Then, for any Borel function 
$F: \Omega_+ \times \Omega_+ \rightarrow \R_+$,
\[
\begin{split}
& \int_{\Gamma(f)} 
\mathrm{d}s \otimes \mathrm{d}a 
\, \frac{1}{\bar s(f,s,a) - \underline{s}(f,s,a)} 
F(\hat f^{s,a}, \check f^{s,a}) \\
& \quad = 
\int_{\Gamma^*(e,U(f))}
\mathrm{d}s \otimes \mathrm{d}a 
\, \frac{1}{\bar s(e,s,a) - \underline{s}(e,s,a)} 
F(\hat e^{s,a}, K^\leftarrow(\cdot; \check e^{s,a},  \check U(e, U(f), s, a)).
\end{split}
\]
\end{corollary}\smallskip

\section{Standard Brownian excursion and length measure}

We first recall a result (Proposition~\ref{P:disintegration} below)
that appears as Corollary 5.2 in \cite{MR2243874}.
It says that if we pick an excursion
$e$ according to the standard excursion distribution
$\PP_{++}$ and then pick  a point $(s,a) \in \Gamma(e)$
according to the $\sigma$-finite measure
\[
\mathrm{d}s \otimes \mathrm{d}a \, \frac{1}{\bar s(e,s,a) - \underline{s}(e,s,a)}
\]
so that the time point $\underline{s}(e,s,a)$ is picked according to 
the $\sigma$-finite measure $m_e$, then the
following objects are independent:
\begin{itemize}
\item[(a)]
the length of the excursion  above level $a$
that straddles time $s$;
\item[(b)]
the excursion obtained by taking the excursion
above level $a$
that straddles time $s$, turning it (by a shift of axes)
into an excursion $\hat e^{s,a}$ above level zero
starting at time zero, 
and then Brownian re-scaling $\hat e^{s,a}$ to
produce an excursion of unit length;
\item[(c)]
the excursion obtained by taking
the excursion $\check e^{s,a}$ that comes from excising
$\hat e^{s,a}$ and closing up the gap, and then
Brownian re-scaling $\check e^{s,a}$
to produce an excursion of unit length;
\item[(d)]
the starting time $\underline{s}(e,s,a)$ of the excursion
above level $a$ that straddles time $s$ rescaled by
the length of  $\check e^{s,a}$ to give
a time in the interval $[0,1]$.
\end{itemize}
Moreover, 
\begin{itemize}
\item
the length in (a) is ``distributed'' according to
the $\sigma$-finite measure 
\[
 \frac{1}{2 \sqrt{2 \pi}}
   \frac{\mathrm{d}r}{\sqrt{(1-r)r^3}}, \quad r \in [0,1];
\]
\item
the unit length excursions in (b) and (c)
are both distributed as standard Brownian excursions
(that is, according to $\PP_{++}$);
\item
the time
in (d) is uniformly distributed on the interval $[0,1]$. 
\end{itemize}

\begin{definition}
For $c > 0$, let $\CS_c: \Omega_{+}^1 \rightarrow \Omega_{+}^c$ be the 
Brownian re-scaling map defined by
\begin{equation*}
\CS_c f := \sqrt{c} f(\cdot/c).
\end{equation*}
\end{definition}\smallskip

\begin{proposition}
\label{P:disintegration}
For any Borel function 
$F: [0,1] \times \Omega_{++} \times \Omega_{++} \rightarrow \R_+$, 
\begin{equation*}
\begin{aligned}
   &\int\PP_{++} (\mathrm{d}e)\,
   \int_{\Gamma(e)}
   \frac{\mathrm{d}s\otimes\mathrm{d}a}{\bar s(e,s,a)-\underline{s}(e,s,a)}
   F\Bigl(\frac{\underline{s}(e,s,a)}{\zeta(\check e^{s,a})}, \hat e^{s,a},\check e^{s,a}\Bigr) 
  \\
   &=
   \int_{[0,1]} \mathrm{d}v \,
   \frac{1}{2 \sqrt{2 \pi}}\int_{[0,1]} 
   \frac{\mathrm{d}r}{\sqrt{(1-r)r^3}}
   \int\PP_{++}(\mathrm{d}e') \otimes \PP_{++}(\mathrm{d}e'')
   \,F(v, \CS_r e', \CS_{1-r} e'').
\end{aligned}
\end{equation*}
\end{proposition}\smallskip

With Proposition~\ref{P:Bertoin_Pitman} and Corollary~\ref{C:length_bridge_length_excursion} in mind, 
we want to obtain an analogous result with $\Gamma(e)$ replaced by $\Gamma^*(e,u)$, where $u$
is picked uniformly from $[0,1]$.

\begin{corollary}
\label{C:disint_with_uniform}
For any Borel function 
$G: [0,1] \times [0,1] \times \Omega_{++} \times \Omega_{++} \rightarrow \R_+$, 
\begin{equation*}
\begin{aligned}
   &\int_{[0,1]} \mathrm{d}u \, \int \PP_{++} (\mathrm{d}e)\,
   \int_{\Gamma^*(e,u)}
   \frac{\mathrm{d}s\otimes\mathrm{d}a}{\bar s(e,s,a)-\underline{s}(e,s,a)} \\
   & \quad \times G\left(\frac{\check U(e,u,s,a)}{\zeta(\check e^{s,a})}, \frac{\underline{s}(e,s,a)}{\zeta(\check e^{s,a})}, \hat e^{s,a},\check e^{s,a}\right) 
  \\
   &=
   \int_{[0,1]} \mathrm{d}u \,
   \int_{[0,1]} \mathrm{d}v \,
   \frac{1}{2 \sqrt{2 \pi}}\int_{[0,1]} 
   \mathrm{d}r \, \sqrt{\frac{1-r}{r^3}} \, 
   \int\PP_{++}(\mathrm{d}e') \otimes \PP_{++}(\mathrm{d}e'') \\
   & \quad \times G(u, v, \CS_r e', \CS_{1-r} e'').
\end{aligned}
\end{equation*}
\end{corollary}\smallskip

\begin{proof}
For $v,r \in [0,1]$ and $u \in [0, (1-r)v[ \, \cup \, ](1-r)v + r , 1]$, put
\[
\breve U(u,v,r)
=
\begin{cases}
\frac{u}{1-r},& 0 \le u < (1-r)v, \\
\frac{u - r}{1-r},&  (1-r)v + r < u \le 1.
\end{cases}
\]

From Proposition~\ref{P:disintegration}, we have
\[
\begin{split}
& \int \PP_{++} (\mathrm{d}e)\,
   \int_{\Gamma^*(e,u)}
   \frac{\mathrm{d}s\otimes\mathrm{d}a}{\bar s(e,s,a)-\underline{s}(e,s,a)} \\
   & \quad \times
   G\left(\frac{\check U(e,u,s,a)}{\zeta(\check e^{s,a})}, \frac{\underline{s}(e,s,a)}{\zeta(\check e^{s,a})}, \hat e^{s,a},\check e^{s,a}\right)\\
& = \int \PP_{++} (\mathrm{d}e)\,
   \int_{\Gamma(e)}
   \frac{\mathrm{d}s\otimes\mathrm{d}a}{\bar s(e,s,a)-\underline{s}(e,s,a)}
   \, \mathbf{1}\{u \notin [\underline{s}(e,s,a), \bar{s}(e,s,a)]\} \\
   & \quad \times G\left(\frac{\check U(e,u,s,a)}{\zeta(\check e^{s,a})}, \frac{\underline{s}(e,s,a)}{\zeta(\check e^{s,a})}, \hat e^{s,a},\check e^{s,a}\right)\\
& =
   \int_{[0,1]} \mathrm{d}v \,
   \frac{1}{2 \sqrt{2 \pi}}\int_{[0,1]} 
   \frac{\mathrm{d}r}{\sqrt{(1-r)r^3}}
   \int\PP_{++}(\mathrm{d}e') \otimes \PP_{++}(\mathrm{d}e'') \\
   & \quad \times \mathbf{1}\{u \notin [(1-r)v, (1-r)v + r]\} \\
& \quad \times G\left(\breve U(u,v,r),  v, \CS_r e', \CS_{1-r} e''\right).\\
\end{split}
\]
The change of variable $w = \breve U(u,v,r)$ gives
\[
\begin{split}
& \int_{[0,1]} \mathrm{d}u \,
\mathbf{1}\{u \notin [(1-r)v, (1-r)v + r]\} 
G\left(\breve U(u,v,r),  v, \CS_r e', \CS_{1-r} e''\right) \\
& \quad =
(1-r)
\int_{[0,1]} \mathrm{d}w \,
G\left(w,  v, \CS_r e', \CS_{1-r} e''\right),\\
\end{split}
\]
and the result follows.
\end{proof}\smallskip

\section{A symmetric measure on $\Omega_+^1 \times \Omega_+^1$}

\begin{definition}
Fix a function $f \in \Omega_+^1$ and suppose that $v \in G(f,a)$
is the starting point of an excursion of $f$ above some level
$a$.  Write
\[
\delta(f,v) := \inf\{t>v: f(t) = a\}
\]
for the time at which the excursion finishes.  Thus,
$\underline{s}(f,s,a) = v$ and $\bar s(f,s,a) = \delta(f,v)$
for any $s \in ]v, \delta(f,v)[$.
Define $\tilde e^{f,v} \in \Omega_{++}$ by
\[
\tilde e^{f,v} := 
\begin{cases}
f(t + v) - f(v),& 0 \le t \le v - \delta(f,v), \\
0,&  t > v - \delta(f,v).
\end{cases}
\]
That is, $\tilde e^{f,v}$ is the result 
of taking the excursion starting and ending
at times $v$ and $\delta(f,v)$, respectively, and shifting the
time and space axes to obtain an excursion that starts at position $0$ at time $0$.
Given $w \in [0,1] \setminus [v, \delta(f,v)]$, denote by
$f^{v,w} \in \Omega_+^1$ the path defined as follows.
If $w>v$ (so that $w > \delta(f,v)$), then
\[
f^{v,w}(t) := 
\begin{cases}
f(t),&  0 \le t < v, \\
f(t - v + \delta(f,v)),&  v \le t < v - \delta(f,v) + w, \\
\tilde e^{f,v}(t - (v - \delta(f,v) + w)) + f(w),&  v - \delta(f,v) + w \le t < w, \\
f(t),&  t \ge w.
\end{cases}
\]
If $w<v$, then
\[
f^{v,w}(t) :=
\begin{cases}
f(t),&  0 \le t < w, \\
\tilde e^{f,v}(t-w) + f(w),&  w \le t < w - v + \delta(f,v), \\
f(t + v - \delta(f,v)),&  w - v + \delta(f,v) \le t < \delta(f,v), \\
f(t),&  t \ge \delta(f,v).\\
\end{cases}
\]
In other words, the excursion of $f$ starting at time $v$ is first moved so that
it starts at $w$ and then the resulting gap left between times $v$ and $\delta(f,v)$
is closed up -- see Figure~\ref{fig:f_to_fvw}.
\begin{figure}[htbp]
	\begin{center}
		\includegraphics[width=1.00\textwidth]{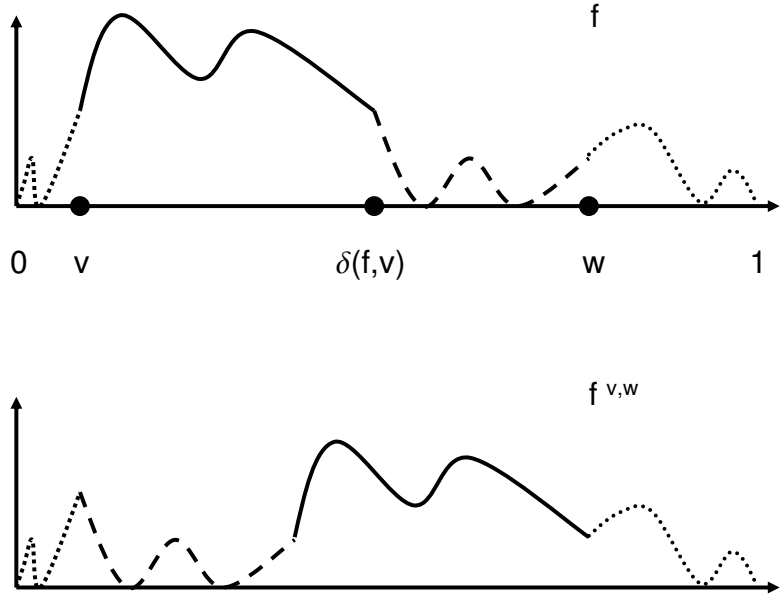}
	\end{center}
	\caption{The transformation taking the path $f$ to the path $f^{v,w}$
	when $w>v$.  The figure for $w<v$ is similar.}
	\label{fig:f_to_fvw}
\end{figure}
\end{definition}\smallskip

\begin{definition}
Define a kernel $\kappa_+$ on $\Omega_+^1$ by
\[
\kappa_+(f,B) := 
\int_0^\infty \mathrm{d}a \, \sum_{v \in {\mathcal G}(f,a)} 
\frac{1}{1 - (\delta(f,v) - v)} \int_{[0,1] \setminus [v,\delta(f,v)]} \mathrm{d}w
\, \mathbf{1}(f^{v,w} \in B).
\]
That is, a starting point $v$ of an excursion is chosen according to
the measure $m_f$ corresponding to length measure $\mu_{T_f}$ on the $\R$-tree associated with
$f$,  this excursion is then relocated so that it starts at a uniformly
chosen point  $w \in [0,1] \setminus [v,\delta(f,v)]$, and finally the
resulting gap is closed up.  
Define a measure $\J_+$ on $\Omega_+^1 \times \Omega_+^1$ by 
\[
\J_+(\mathrm{d}f', \mathrm{d} f'') := \PP_+(\mathrm{d}f') \kappa_+(f', \mathrm{d}f'').
\]
\end{definition}\smallskip

\begin{proposition}
\label{P:symmetry_J+}
The measure $\J_+$ is symmetric.
\end{proposition}\smallskip

\begin{proof} Given $e', e'' \in \Omega_{++}^1$, $v \in [0,1]$,
and $r \in ]0, 1]$, define 
$e^\circ(\cdot; e', e'', v, r) \in \Omega_{++}^1$ by
\[
\begin{split}
& e^\circ(t; e', e'', v, r) \\
& \quad :=
\begin{cases}
\CS_{1-r}e''(t), &
   0 \le t \le (1-r)v, \\
\CS_{1-r} e''((1-r)v) + \CS_r e'(t - (1-r)v), &
   (1-r)v \le t \le (1-r)v + r, \\
\CS_{1-r} e''(t - r), & 
   (1-r)v + r \le t \le 1. 
\end{cases}
\\
\end{split}
\]
That is, $e^\circ(\cdot; e', e'', v, r)$
is the excursion that arises from Brownian re-scaling
$e'$ and $e''$ to have lengths $r$ and $1-r$,
respectively, and then inserting the re-scaled version
of $e'$ into the re-scaled version of $e''$ at a position
that is a fraction $v$ of the total length of the re-scaled
version of $e''$.

Also, for $u \in [0,1]$ set
\[
\tilde U(u,v,r)
:=
\begin{cases}
(1-r)u,& 0 \le u \le v, \\
r + (1-r)u,& v < u \le 1,
\end{cases}
\]
so that $\tilde U(u,v,r)$ belongs to the
set $[0,(1-r)v[ \, \cup \, ](1-r)v + r, 1]$
for Lebesgue almost all $u \in [0,1]$
and the push-forward of Lebesgue measure on $[0,1]$
by the map $u \mapsto \tilde U(u,v,r)$ is  the
uniform distribution on this union of two intervals.

Define a measure $\J_{++}$ on $[0,1] \times [0,1] \times \Omega_{++}^1 \times \Omega_{++}^1$ by
\[
\begin{split}
& \int 
\J_{++}(\mathrm{d} u^*, \mathrm{d} u^{**}, \mathrm{d} e^*, \mathrm{d} e^{**})
G(u^*, u^{**}, e^*, e^{**}) \\
& \quad :=
\int_{[0,1]^3} \mathrm{d}u \otimes \mathrm{d}v \otimes \mathrm{d}w
\frac{1}{2 \sqrt{2 \pi}} 
\int_{[0,1]} \mathrm{d} r \sqrt{\frac{1-r}{r^3}}
\int \PP_{++}(\mathrm{d}e') \otimes \PP_{++}(\mathrm{d}e'') \\
& \qquad \times G\left(\tilde U(u,v,r), \tilde U(u,w,r), e^\circ(\cdot; e', e'', v, r), 
e^\circ(\cdot; e', e'', w, r)\right) \\
\end{split}
\]
for any non-negative Borel function $G$.

Clearly, the measure $\J_{++}$ is preserved by pushing it forward with the map
$(u^*, u^{**}, e^*, e^{**}) \mapsto (u^{**}, u^*,  e^{**}, e^*)$. Also, it follows
from Lemma~\ref{L:K_arrow_relations},  Proposition~\ref{P:Bertoin_Pitman},
Corollary~\ref{C:length_bridge_length_excursion} and Corollary~\ref{C:disint_with_uniform}
that the measure $\J_+$ is the push-forward of the measure $\J_{++}$ by the map
\[
(u^*, u^{**}, e^*, e^{**}) 
\mapsto 
(K^\leftarrow(\cdot; e^*,  u^*), K^\leftarrow(\cdot; e^{**},  u^{**})),
\]
and the result follows.
\end{proof}\smallskip

By construction, the measure $\J_+$ is concentrated on pairs $(f', f'')
\in \Omega_+^1 \times \Omega_+^1$ such that $f''$ is obtained from $f'$
by the re-location of an excursion.  If we shift the starting point
of this excursion in space and time to the origin to obtain
an element of $\Omega_{++}$, then the $\sigma$-finite law 
of this shifted excursion is
\[
\begin{split}
& \PP_+
\left[
\int_0^\infty \mathrm{d}a \, \sum_{v \in {\mathcal G}(f,a)}
\, \mathbf{1}\left(\tilde e^{f,v} \in \cdot\right) 
\right] \\
& \quad =
\PP_+
\left[
\int_{\Gamma(f)} 
\mathrm{d}s \otimes \mathrm{d}a 
\, \frac{1}{\bar s(f,s,a) - \underline{s}(f,s,a)} 
\, \mathbf{1}\left(\hat f^{s,a} \in \cdot\right)
\right]. \\
\end{split}
\]
Informally, this is the law of the excursion under $\J_+$, but
we note while $\J_+$ is concentrated on pairs $(f',f'')$
of the form $(f,f^{v,w})$ for some $v,w \in [0,1]$, the value
of $v$ and  the corresponding excursion $\tilde e^{f,v}$ cannot
be uniquely reconstructed from $(f',f'')$.
Arguing as in the proof of Proposition \ref{P:symmetry_J+}, this
law is given by
\[
\frac{1}{2 \sqrt{2 \pi}} 
\int_{[0,1]} \mathrm{d} r \, \sqrt{\frac{1-r}{r^3}}
\PP_{++}
\left\{e \in \Omega_{++}^1 :
\mathcal{S}_r(e) \in \cdot
\right\}.
\]
We need the following properties of this law.

\begin{proposition}
\label{P:shifted_excursion}
\begin{itemize}
\item[(i)]
For $0<t\le 1$,
\[
\int \PP_+({\mathrm d}f)
\int_0^\infty \mathrm{d}a \, \sum_{v \in {\mathcal G}(f,a)}
\, \mathbf{1}(\zeta(\tilde e^{f,v}) > t) 
 =
\frac{1}{\sqrt{2 \pi }}
\left(
\sqrt{\frac{1}{t}-1} + \arcsin\left(\sqrt{t}\right) - \frac{\pi}{2}
\right),
\]
and hence
\[
\int \PP_+({\mathrm d}f)
\int_0^\infty \mathrm{d}a \, \sum_{v \in {\mathcal G}(f,a)}
(\zeta(\tilde e^{f,v}))^2 
=\frac
{\pi^{\frac{1}{2}}}
{16 \sqrt{2}}.
\]
\item[(ii)]
For $x>0$,
\[
\int \PP_+({\mathrm d}f)
\int_0^\infty \mathrm{d}a \, \sum_{v \in {\mathcal G}(f,a)}
\, \mathbf{1}\left(\max(\tilde e^{f,v}) > x\right) 
=
\sum_{n=1}^\infty
\int_{2 n x}^\infty \mathrm{d}z \, \exp\left(-\frac{z^2}{2}\right),
\]
and hence
\[
\int \PP_+({\mathrm d}f)
\int_0^\infty \mathrm{d}a \, \sum_{v \in {\mathcal G}(f,a)}
\left(\max(\tilde e^{f,v})\right)^2 
= 
\frac{\pi^{\frac{5}{2}}}{24 \sqrt{2}}.
\]
\end{itemize}
\end{proposition}

\begin{proof}
(i) By the remarks prior the statement of the proposition, the
quantity in the first claim is
\[
\frac{1}{2 \sqrt{2 \pi}} 
\int_{[t,1]} \mathrm{d} r \, \sqrt{\frac{1-r}{r^3}},
\]
and a straightforward integration shows that this has
the stated value.  The second claim follows by an equally
straightforward integration by parts.

\noindent
(ii) Again by the remarks prior the statement of the proposition, the
quantity in the first claim is
\[
\frac{1}{2 \sqrt{2 \pi}} 
\int_{[0,1]} \mathrm{d} r \, \sqrt{\frac{1-r}{r^3}}
\PP_{++}\left\{e \in \Omega_{++}^1 : \max(e) > \frac{x}{\sqrt{r}}\right\}.
\]
From Theorem 5.2.10 in \cite{MR82m:60098}, we have that
\[
\PP_{++}\left\{e \in \Omega_{++}^1 : \max(e) > y \right\}
= 2 \sum_{n=1}^\infty (4 n^2 y^2 - 1) \exp(-2 n^2 y^2),
\]
and an integration establishes the claim.

An integration by parts shows that the quantity in the second claim is
\[
\sqrt{\frac{\pi}{32}} \sum_{n=1}^\infty \frac{1}{n^2} 
= \sqrt{\frac{\pi}{32}} \frac{\pi^2}{6},
\]
as required.
\end{proof}

\section{A Dirichlet form}
\label{S:Dirichlet}

Recall that any $f \in \Omega_+^1$ is associated with a
$\R$-tree $(T_f, d_{T_f})$ that arises
as a quotient of $[0,1]$ under an equivalence
relation defined by $f$.  Moreover, we may 
equip this $\R$-tree with the root
$\rho_{T_f}$ that is the image of $0$ 
under the quotient and the weight $\nu_{T_f}$
that is the push-forward of Lebesgue measure on $[0,1]$ by the
quotient map.

\begin{definition}
Define the probability measure ${\mathbf P}$
on ${\mathbf T}^{\mathrm{wr}}$ to be the push-forward
of the probability measure $\PP_+$ on 
$\Omega_+^1$ by the map 
\[
f \mapsto (T_{2f}, d_{T_{2f}}, \rho_{T_{2f}}, \nu_{T_{2f}}).
\]
Define the measure ${\mathbf J}$ on 
${\mathbf T}^{\mathrm{wr}} \times {\mathbf T}^{\mathrm{wr}}$
to be {\bf twice} the push-forward of $\J_+$ by the map
\[
(f',f'') \mapsto
( 
(T_{2f'}, d_{T_{2f'}}, \rho_{T_{2f'}}, \nu_{T_{2f'}})
,
(T_{2f''}, d_{T_{2f''}}, \rho_{T_{2f''}}, \nu_{T_{2f''}})
).
\]
\end{definition}\smallskip

\begin{proposition}
\label{P:properties_jump}
\begin{itemize}
\item[(i)] The measure ${\mathbf J}$ is symmetric.
\item[(ii)] For each compact subset $K \subset{\mathbf T}^{\mathrm{wr}}$ 
and open
subset $U$ such that 
$K \subset U \subseteq{\mathbf T}^{\mathrm{wr}}$,
\begin{equation*}
   {\mathbf J}(K \times ({\mathbf T}^{\mathrm{wr}}\setminus U))
 <
   \infty.
\end{equation*}
\item[(iii)] The function $\Delta_{\mathrm{GH^{wr}}}$ is square-integrable
with respect to $J$, that is,
\begin{equation*}
   \int
   {\mathbf J}(\mathrm{d}T',\mathrm{d}T'') \,
   \Delta_{\mathrm{GH^{wr}}}^2(T',T'')
   <\infty.
\end{equation*}
\end{itemize}
\end{proposition}\smallskip

\begin{proof}
(i) This is immediate from Proposition~\ref{P:symmetry_J+}.

\noindent
(ii) By construction, the measure ${\mathbf J}$
has the following description.  Firstly, a weighted
rooted compact $\R$-tree $T' \in {\mathbf T}^{\mathrm{wr}}$
is chosen according to ${\mathbf P}$.  A point $v \in T'$
is chosen according to the length measure $\mu_{T'}$ and
another point $w \in T'$ is chosen according to the
renormalization of the weight $\nu_{T'}$ outside of the
subtree $S^{T',v}$ of points ``above'' $v$
(that is, of points $x$ such that $v$ belongs to the
segment $[\rho_{T'}, x[$).  The subtree $S^{T',v}$
is then pruned off and re-attached at $w$ to form
a new $\R$-tree $T''$.  More formally, the $\R$-tree
$T''$ can be identified as the set $T'$ equipped with
new metric $d_{T''}$ given by
\[
   d_{T''}(x,y)
 :=
 \begin{cases}
 d(x,y),& x,y\in S^{T',v},\\
   d(x,y), & x,y\in T'\setminus S^{T',v},\\
   d(x,v)+d(w,y),& x\in S^{T',v}, \, y\in T' \setminus
   S^{T',v},\\
   d(y,v)+d(w,x),& y\in S^{T',v}, x\in T'\setminus S^{T',v}.
 \end{cases}
\]
With this identification, $\rho_{T''} = \rho_{T'}$ and
$\nu_{T''} = \nu_{T'}$.

We claim that if, for some $\varepsilon>0$, 
\[
\max_{x \in S^{T',v}} d_{T'}(v,x) \le \varepsilon
\]
and
\[
\nu_{T'}(S^{T',v}) \le \varepsilon,
\]
 then
$\Delta_{\mathrm{GH^{wr}}}(T',T'') \le \varepsilon$.
Firstly, the map $f:T' \rightarrow T''$ defined by 
\[
f(x) :=
\begin{cases} 
x,& x \in T' \setminus S^{T',v},\\
u, & x \in S^{T',v},
\end{cases}
\]
is such that $f(\rho_{T'}) = \rho_{T''}$ and
\[
\sup\{|d_{T'}(x,y) - d_{T''}(f(x),f(y))| : x,y \in T'\} \le \varepsilon.
\]
Moreover, it is immediate that
\[
d_P(f_*\nu_{T'}, \nu_{T''}) \le \varepsilon,
\]
and so $f \in F_{T',T''}^\varepsilon$.  Note also that
$S^{T'',w} = S^{T',v}$ as sets,
\[
\max_{x \in S^{T'',w}} d_{T''}(w,x) = \max_{x \in S^{T',v}} d_{T'}(v,x),
\]
and
\[
\nu_{T''}(S^{T'',w}) = \nu_{T'}(S^{T',v}),
\]
and so a similar argument shows that the map
$g:T'' \rightarrow T'$ defined by 
\[
g(x) :=
\begin{cases} 
x,& x \in T'' \setminus S^{T'',w},\\
w, & x \in S^{T'',w}
\end{cases}
\]
belongs to $F_{T'',T'}^\varepsilon$.
Thus, $\Delta_{\mathrm{GH^{wr}}}(T',T'') \le \varepsilon$
as required.

Now let $K$ and $U$ be as in the statement of part (ii).
The result is trivial if $K = \emptyset$,
so we assume that $K \ne \emptyset$.
Since ${\mathbf T}^{\mathrm{wr}}\setminus U$ and 
$K$ are disjoint closed sets and $K$ 
is compact, we have that
\[
   c
 :=
   \inf_{T'\in K , T''\in U} \Delta_{\mathrm{GH^{wr}}}(T',T'')>0.
\]
From what we have just observed, if $0 < \varepsilon < c$, then,
by Proposition~\ref{P:shifted_excursion},
\[
\begin{split}
   {\mathbf J}(K \times ({\mathbf T}^{\mathrm{wr}}\setminus U) )
& \le
{\mathbf J}\left\{(T',T'') : \Delta_{\mathrm{GH^{wr}}}(T',T'') > \varepsilon\right\} \\
& \le
\int {\mathbf P}({\mathrm d}T') 
\int_{T'} \mu_{T'} ({\mathrm d}v)
\, {\mathbf 1}\left(\max_{x \in S^{T',v}} d_{T'}(v,x) > \varepsilon \right) \\
& \quad +
\int {\mathbf P}({\mathrm d}T') 
\int_{T'} \mu_{T'} ({\mathrm d}v)
\, {\mathbf 1}\left(\nu(S^{T',v}) > \varepsilon \right) \\
& =
\int \PP_+({\mathrm d}f) \,
2 \int_0^\infty \mathrm{d}a \, \sum_{v \in {\mathcal G}(f,a)}
\, \mathbf{1}\left(\max(\tilde e^{f,v}) > \varepsilon / 2\right) \\
& \quad +
\int \PP_+({\mathrm d}f) \,
2 \int_0^\infty \mathrm{d}a \, \sum_{v \in {\mathcal G}(f,a)}
\, \mathbf{1}\left(\zeta\left(\tilde e^{f,v}\right) > \varepsilon\right) \\
& < \infty, \\
\end{split}
\]
as required.

\noindent
(iii) By an argument similar to that in part (ii) and
Proposition~\ref{P:shifted_excursion}, we have
\[
\begin{split}
& \int
   {\mathbf J}(\mathrm{d}T',\mathrm{d}T'') \,
   \Delta_{\mathrm{GH^{wr}}}^2(T',T'') \\
& \quad = 
\int_0^\infty {\mathrm d}\varepsilon \, 2 \varepsilon \,
{\mathbf J}\left\{(T',T'') : \Delta_{\mathrm{GH^{wr}}}(T',T'') > \varepsilon\right\} \\
& \quad \le
\int_0^\infty {\mathrm d}\varepsilon \, 2 \varepsilon
\int {\mathbf P}({\mathrm d}T') \,
\int_{T'} \mu_{T'} ({\mathrm d}v) \,
{\mathbf 1}\left(\max_{x \in S^{T',v}} d_{T'}(v,x) > \varepsilon \right) \\
& \qquad +
\int_0^\infty {\mathrm d}\varepsilon \, 2 \varepsilon
\int {\mathbf P}({\mathrm d}T') \, 
\int_T' \mu_{T'} ({\mathrm d}v) \,
{\mathbf 1}\left(\nu(S^{T',v}) > \varepsilon \right) \\
& \quad =
\int_0^\infty {\mathrm d}\varepsilon \, 2 \varepsilon
\int \PP_+({\mathrm d}f) \, 
2 \int_0^\infty \mathrm{d}a \, \sum_{v \in {\mathcal G}(f,a)}
\, \mathbf{1}\left(\max(\tilde e^{f,v}) > \varepsilon / 2\right) \\
& \qquad +
\int_0^\infty {\mathrm d}\varepsilon \, 2 \varepsilon
\int \PP_+({\mathrm d}f) \, 
2 \int_0^\infty \mathrm{d}a \, \sum_{v \in {\mathcal G}(f,a)}
\, \mathbf{1}(\zeta(\tilde e^{f,v}) > \varepsilon) \\
& \quad =
\int \PP_+({\mathrm d}f) \,
8 \int_0^\infty \mathrm{d}a \, \sum_{v \in {\mathcal G}(f,a)}
\left(\max(\tilde e^{f,v})\right)^2 \\
& \qquad +
\int \PP_+({\mathrm d}f) \,
2 \int_0^\infty \mathrm{d}a \, \sum_{v \in {\mathcal G}(f,a)}
(\zeta(\tilde e^{f,v}))^2 \\
& \quad < \infty,\\
\end{split}
\]
as required.
\end{proof}

\begin{definition}
Define a bilinear form 
\[
   {\mathcal E}(f,g)
 :=
   \int
  {\mathbf J}(\mathrm{d}T',\mathrm{d}T'')
  \left(f(T'')-f(T')\right)\left(g(T'')-g(T')\right),
\]
for $f,g$ in the domain 
\[
{\mathcal D}^*({\mathcal E})
:=
   \{f\in L^2({\mathbf T}^{\mathrm{wr}},{\mathbf P}):\,\text{$f$ is
     Borel and ${\mathcal E}(f,f)<\infty$}\}.
\]
Here, as usual, $L^2({\mathbf T}^{\mathrm{wr}}, {\mathbf P})$ is equipped with the inner product 
\[
(f,g)_{\mathbf P}:=\int{\mathbf P}(\mathrm{d}x) \, f(x)g(x).
\]
\end{definition}\smallskip

\begin{definition}
Let $\mathcal L$ denote the collection of functions
$f: {\bf T}^{\mathrm{wr}} \rightarrow \R$
such that 
\[
\sup_{T \in {\bf T}^{\mathrm{wr}}}
|f(T)| < \infty
\]
and
\[
\sup_{T', T'' \in {\bf T}^{\mathrm{wr}}, \, T' \ne T''}
\frac{|f(T') - f(T'')|} 
{\Delta_{\mathrm{GH^{wr}}}(T',T'')}
< \infty.
\]
\end{definition}\smallskip

Part (i) (respectively, parts (ii) and (iii))
of the following result may be proved in the same manner
as Lemma~7.1 (respectively, Lemma~7.2 and Theorem~7.3)
of \cite{MR2243874}, with Proposition~\ref{P:properties_jump}
above playing the role played in  \cite{MR2243874}
by Lemma~6.2 of that paper.  Our setting is slightly different,
in that we are working with weighted rooted compact $\R$-trees
rather than just weighted compact $\R$-trees, but this
doesn't require any significant changes in the arguments.
In particular, the analogues of Lemmas~7.5, 7.6 and 7.7
of \cite{MR2243874} go through quite straightforwardly to
establish the tightness property required for 
part (iii) to hold.

\begin{theorem}
\label{T:main}
\begin{itemize}
\item[(i)]
The form 
$({\mathcal E},{\mathcal D}^*({\mathcal E}))$ 
is Dirichlet (that is, it is symmetric, non-negative definite,
Markovian, and closed).
\item[(ii)] 
The set $\mathcal L$ is a vector lattice and an algebra,
and ${\mathcal L} \subseteq {\mathcal D}^*({\mathcal E})$.
Hence, if ${\mathcal D}({\mathcal E})$ denotes the
closure of $\mathcal L$ in ${\mathcal D}^*({\mathcal E}))$,
then $({\mathcal E}, {\mathcal D}({\mathcal E}))$
is also a Dirichlet form.
\item[(iii)]
There is a recurrent ${\mathbf P}$-symmetric Hunt process
$X=(X_t,{\mathbb P}^T)$ on ${\mathbf T}^{\mathrm{wr}}$ 
with Dirichlet form  
$({\mathcal E}, {\mathcal D}({\mathcal E}))$.
\end{itemize}
\end{theorem}\smallskip

\providecommand{\bysame}{\leavevmode\hbox to3em{\hrulefill}\thinspace}
\providecommand{\MR}{\relax\ifhmode\unskip\space\fi MR }
\providecommand{\MRhref}[2]{%
  \href{http://www.ams.org/mathscinet-getitem?mr=#1}{#2}
}
\providecommand{\href}[2]{#2}

\end{document}